\documentclass{article}
\usepackage{amssymb,amsmath,amsthm,graphicx}
\usepackage{xcolor}
\usepackage[all,color]{xy}

\def\qed{\hfill {\hbox{${\vcenter{\vbox{               
   \hrule height 0.4pt\hbox{\vrule width 0.4pt height 6pt
   \kern5pt\vrule width 0.4pt}\hrule height 0.4pt}}}$}}}

\def\tr{\triangleright}
\def\otr{\, \overline{\tr}\, }
\def\utr{\, \underline{\tr}\, }

\newtheorem{theorem}{Theorem}

\newtheorem{proposition}[theorem]{Proposition}

\theoremstyle{definition}
\newtheorem{example}{Example}
\newtheorem{definition}{Definition}
\newtheorem{remark}{Remark}

\date{}

\title{\Large \textbf{Virtual Biquandle Cocycle Quiver Representations}}

\author{Alexander Bishop\footnote{Email: abishop33@students.claremontmckenna.edu.}
\and Jose Ceniceros\footnote{Email: jcenicer@hamilton.edu.}
\and Sam Nelson\footnote{Email: Sam.Nelson@cmc.edu. Partially supported by Simons Foundation collaboration grant 702597.}}

\begin{document}
\maketitle

\begin{abstract}
We introduce quiver representation-valued invariants of oriented virtual knots
and links associated to a choice of finite virtual biquandle, abelian group, 
set of virtual Boltzmann weights, commutative unital ring and set of virtual 
biquandle endomorphisms. As an application we define new infinite families
of polynomial virtual knot and link invariants via decategorification.
\end{abstract}

\parbox{5.5in} {\textsc{Keywords:} Virtual biquandles, Cocycle invariants,
coloring quivers, quiver representations

\smallskip

}

\section{Introduction}\label{I}

\textit{Virtual knots} are knots in thickened surfaces regarded up to 
stabilization of the surface. Classical knot theory can be regarded as 
the special case of virtual knot theory in which the supporting surface 
has genus zero. Since it is impractical to use surfaces of nonzero genus 
for drawing our diagrams, we represent genus in the underlying surface 
with \textit{virtual crossings} which interact with classical crossings 
via the \textit{Virtual Reidemeister moves}. See \cite{K,EN} for more.

\textit{Biquandles} are algebraic structures with axioms arising from the
Reidemeister moves in classical knot theory. Every virtual knot or link has a 
\textit{fundamental biquandle}; an infinite set of computable integer-valued 
knot and link invariants is defined by considering the \textit{homset} 
$\mathrm{Hom}(\mathcal{B}(L),X)$ of biquandle homomorphisms from the 
fundamental biquandle $\mathcal{B}(L)$ of a link $L$ to a finite quandle $X$. 
A choice of diagram for an oriented link $L$, analogously to a choice of 
basis for a vector space, provides a concrete representation of homset 
elements as \textit{biquandle colorings} of the diagram. See \cite{EN} 
for more.

In \cite{CN}, subsets of the set of endomorphisms $\mathrm{Hom}(X,X)$ of 
a finite quandle were used to enhance and categorify the quandle homset 
invariant, yielding directed graph-valued invariants called \textit{quandle 
coloring quivers}. Quivers are small categories, and further invariants 
have been derived from these quivers via decategorification. In \cite{CN2}, 
cocycle values were used to define weighted quivers which categorify the
quandle 2-cocycle invariant. In \cite{SN1} a quiver representation-valued
invariant was defined using quandle coloring quivers and 2-cocycles, with
several families of polynomial invariants defined by decategorification.

In \cite{KM}, the combinatorial diagrammatic interpretation of biquandle 
colorings was used to define \textit{virtual biquandles} for coloring virtual 
knots and links. In \cite{JCN1} cocycles in certain cohomology theories for 
virtual biquandles were used to enhance the counting invariant. In \cite{NT} 
the fundamental virtual kei of unoriented virtual knots were studied. 

In this paper we generalize the biquandle cocycle quiver representation 
from \cite{SN1} to the case of virtual biquandles, categorifying the virtual 
biquandle cocycle invariants from \cite{JCN1} and obtaining new computable 
invariants as a result. The paper is organized as follows. In Section 
\ref{V} we review the notions of virtual knot theory and virtual biquandles. In 
Section \ref{QQ} we define virtual biquandle cocycle quiver representations 
and introduce the new invariants. In Section \ref{EC} we collect some examples 
and computations. We conclude in Section \ref{Q} with some questions for 
future research.

This paper, including all text, diagrams and computational code,
was produced strictly by the authors without the use of generative AI in any 
form.

\section{\textbf{Virtual Biquandles and 2-Cocycles}}\label{V}

We begin with a definition; see \cite{CN,EN} and the references therein for more.

\begin{definition}
A \textit{virtual biquandle} is a set $X$ with two binary operations
$\utr,\otr:X\times X\to X$ and an invertible map $v:X\to X$ satisfying 
the conditions
\begin{itemize}
\item[(i)] For all $x\in X$, $x\utr x=x\otr x$,
\item[(ii)] For all $y\in X$ the maps $\alpha_y,\beta_y:X\to X$ defined by
$\alpha_y(x)=x\otr y$ and $\beta_y(x)=x\utr y$ and the map
$S:X\times X\to X\times X$ defined by $S(x,y)=(\alpha_x(y),\beta_y(x))$
are invertible, 
\item[(iii)] For all $x,y,z\in X$ we have 
\[\begin{array}{rcl}
(x\utr y)\utr (z\utr y) & = & (x\utr z)\utr (y\utr z), \\
(x\utr y)\otr (z\utr y) & = & (x\otr z)\utr (y\otr z)\ \mathrm{and}\\
(x\otr y)\otr (z\otr y) & = & (x\otr z)\otr (y\utr z), 
\end{array}\]
and
\item[(iv)] For all $x,y\in X$ we have 
\[v(x\utr y)=v(x)\utr v(y) \quad \mathrm{and}\quad v(x\otr y)=v(x)\otr v(y)\]
\end{itemize}
A virtual biquandle in which $v=\mathrm{Id}_X$ is called a \textit{biquandle}.
\end{definition}

\begin{example}
Any module over $\mathbb{Z}[t^{\pm 1},s^{\pm 1}]$ is a virtual biquandle with
\[x\utr y=tx+(s-t)y,\quad x\otr y=sx\quad \mathrm{and}\quad v(x)=ax\]
for any unit $a\in X$.
\end{example}

\begin{example}
Any group $G$ is a virtual biquandle with
\[x\utr y=y^{-1}xy^{-1},\quad x\otr y=x^{-1}\quad \mathrm{and}\quad v(x)=x^{-1}.\]
\end{example}

\begin{example}
Every oriented virtual knot or link $L$ has a \textit{fundamental virtual 
biquandle} $\mathcal{V}(L)$ with presentation given by a set of generators 
corresponding to virtual semi-arcs in a virtual diagram representing $L$ with 
relations at crossings as shown:
\[\includegraphics{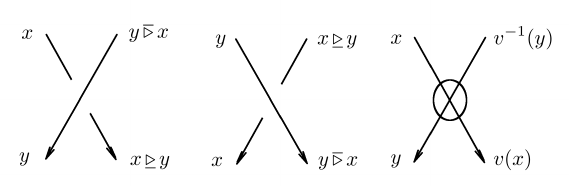}\]
The elements of $\mathcal{V}(L)$ are equivalence classes of virtual biquandle
words in these generators modulo the equivalence relation determined by the
virtual biquandle axioms and the crossing relations. The biquandle axioms
are chosen so that virtual Reidemeister moves on virtual link diagrams
induce Tietze moves on the presentation, and hence $\mathcal{V}(L)$ is
an invariant of oriented virtual knots and links.
\end{example}

\begin{example}\label{ex1}
We can specify a finite virtual biquandle by giving its operation
tables. For example, the following tables describe a virtual biquandle
with three elements:
\[
\begin{array}{r|rrr}
\utr & 1 & 2 & 3 \\ \hline
1 & 2 & 2 & 2 \\
2 & 3 & 3 & 3 \\
3 & 1 & 1 & 1 
\end{array}
\quad
\begin{array}{r|rrr}
\otr & 1 & 2 & 3 \\ \hline
1 & 2 & 1 & 3 \\
2 & 1 & 3 & 2 \\
3 & 3 & 2 & 1 
\end{array}
\quad
\begin{array}{r|c}
x & v(x) \\ \hline
1 & 2 \\
2 & 3  \\
3 & 1 \\
\end{array}.
\]
\end{example}

\begin{definition}
Let $(X,v)$ be a virtual biquandle and $L$ an oriented virtual knot or
link. Then the set of virtual biquandle homomorphisms from $\mathcal{V}(L)$
to $X$, denoted $\mathrm{Hom}(\mathcal{V}(L),X)$, is the \textit{virtual 
biquandle homset} of $L$ with respect to $(X,v)$.
\end{definition}

We then have the following standard result; see \cite{KM}.

\begin{theorem}
If $L$ and $L'$ are equivalent oriented virtual links, then
there exists a bijection between $\mathrm{Hom}(\mathcal{V}(L),(X,v))$
and $\mathrm{Hom}(\mathcal{V}(L'),(X,v))$.
\end{theorem}

More precisely, a \textit{virtual biquandle coloring} or 
\textit{$(X,v)$-coloring} of a virtual knot diagram $D$ is an assignment
of an element of $X$ to each semiarc in $D$ such that at every crossing
we have
\[\includegraphics{ab-jc-sn-1.pdf}.\]
Then $X$-colorings of $D$ uniquely determine and are determined by
elements of $\mathrm{Hom}(\mathcal{V}(L),(X,v))$. In particular, 
$(X,v)$-colored diagrams represent elements of 
$\mathrm{Hom}(\mathcal{V}(L),(X,v))$ analogously to the way matrices represent
linear transformations with respect to a choice of basis, with the role of 
change-of-basis matrices played by $(X,v)$-colored virtual Reidemeister moves
on $D$.

\begin{example}
Let $(X,v)$ be the virtual biquandle with operation tables
\[
\begin{array}{r|rrr}
\utr & 1 & 2 & 3 \\ \hline
1 & 2 & 2 & 2 \\
2 & 3 & 3 & 3 \\
3 & 1 & 1 & 1 
\end{array}
\quad
\begin{array}{r|rrr}
\otr & 1 & 2 & 3 \\ \hline
1 & 2 & 1 & 3 \\
2 & 1 & 3 & 2 \\
3 & 3 & 2 & 1 
\end{array}
\quad
\begin{array}{r|c}
x & v(x) \\ \hline
1 & 2 \\
2 & 3  \\
3 & 1 \\
\end{array}.
\] 
Then the 
virtual knot $2.1$ has homset
\[\mathrm{Hom}(\mathcal{V}(2.1),(X,v))
=\left\{
\raisebox{-0.5in}{\scalebox{0.8}{\includegraphics{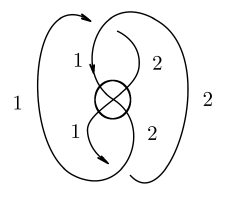}}},\
\raisebox{-0.5in}{\scalebox{0.8}{\includegraphics{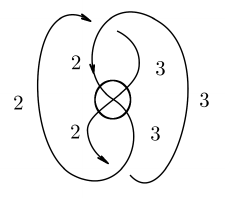}}},\
\raisebox{-0.5in}{\scalebox{0.8}{\includegraphics{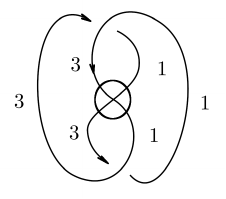}}}
\right\}.\]
\end{example}

\begin{definition}
Let $(X,\utr, \otr,v),(Y, \,\underline{*}\,, \,\overline{*}\,, v')$ be 
virtual biquandles. A \textit{virtual biquandle homomorphism} is a map 
$\sigma: X \rightarrow Y$ such that for all $x,y \in X$ we have
\[ \sigma(x\utr y)=\sigma(x)\,\underline{*}\,\sigma(y),\quad
\sigma(x\otr y)=\sigma(x)\,\overline{*}\,\sigma(y)\quad \mathrm{and}\quad 
v'(\sigma(x))=\sigma(v(x)).\]
We denote the set of virtual biquandle homomorphisms as $\mathrm{Hom}(X,Y)$. 
A bijective virtual biquandle homomorphism $\sigma: X \rightarrow Y$ 
is a \textit{virtual biquandle isomorphism}.
A virtual biquandle homomorphism $\sigma: X \rightarrow X$ is a 
\textit{virtual biquandle endomorphism}. We denote the set of virtual 
biquandle endomorphisms as $\mathrm{End}(X)$.
A bijective virtual biquandle endomorphism $\sigma: X \rightarrow X$ 
is a \textit{virtual biquandle automorphism}. We denote the set of 
virtual biquandle automorphisms as $\mathrm{Aut}(X)$.
\end{definition}

\begin{remark}
A virtual biquandle is a biquandle with a distinguished automorphism $v$.
Then a virtual biquandle endomorphism is a biquandle endomorphism which
commutes with $v$.
\end{remark}

Next, we have a definition from \cite{JCN1}, updated for our current notation:

\begin{definition}
Let $X$ be a virtual biquandle and $A$ an abelian group. A \textit{virtual 
biquandle Boltzmann weight} is given by a pair of functions 
$\phi,\psi:X\times X\to A$ such that for all $x,y,z\in X$ we have
\[
\begin{array}{rcll}
\phi(x,y)+\phi(y,z)+\phi(x\utr y,z\otr y)
-\phi(x,z)-\phi(y\otr x, z\otr x)-\phi(x\utr z, y\utr z) & = & 0 & (i.i)\\
\psi(x,y)+\psi(y,z)+\psi(v(x),v^{-1}(z))-\psi(x,z)
-\psi(v^{-1}(y),v^{-1}(z))-\psi(v(x),v(y)) & = & 0 & (i.ii)\\
\phi(x,y)+\psi(y,z)+\psi(x\utr y, v^{-1}(z))-\psi(x,z)-\phi(v(x),v(y))
-\psi(y\otr x,v^{-1}(z)) & = & 0 & (i.iii)\\
\end{array}
\]
and for all $x\in X$ we have
\[\phi(x,x)=0,\quad \psi(x,x)=0 \quad \mathrm{and}\quad \psi(x,y)+\psi(y,x)=0.\]
If in addition, for all $x,y\in X$ we have 
\[\phi(x,y)=\phi(v(x),v(y))\]
then we say $\phi,\psi$ are \textit{strongly compatible}.
\end{definition}

The motivation for this definition is the following. 
Let $(X,v)$ be a virtual biquandle and $L$ an oriented virtual knot
or link represented by a diagram $D$. Let $A$ be an abelian group and 
$\phi,\psi$ be a virtual biquandle Boltzmann weight taking values in $A$.
The virtual biquandle Boltzmann weight definition is chosen so that the
sum of contributions from $(X,v)$-colored classical and virtual crossings
\[\includegraphics{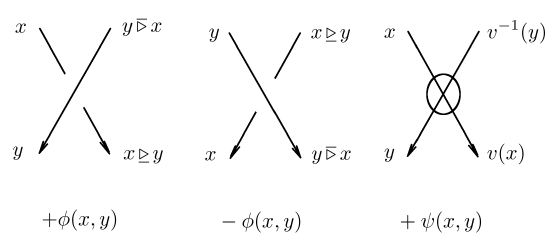}\]
is unchanged by virtual Reidemeister moves. In the case of strong 
compatibility, the classical ($\phi$) and virtual ($\psi$) parts of the 
Boltzmann weight are independently invariant, so the pair of these weights 
is invariant.

We then have the following theorem; see also \cite{JCN1}.

\begin{theorem}
The multiset of Boltzmann weight values over the virtual biquandle homset
is an invariant of oriented classical and virtual knots and links. In the 
case of strong compatibility, the multiset of pairs of classical and virtual 
Boltzmann weights over the virtual biquandle homset is an invariant of oriented
classical and virtual knots and links.
\end{theorem}

We can specify a Boltzmann weight as a pair of tables whose entries are the
coefficients for $\chi_{x,y}$ in $\phi$ and $\psi$.

\begin{example}
Let $X$ be the virtual biquandle specified by the operation tables
\[
\begin{array}{r|rrr}
\utr & 1 & 2 & 3 \\ \hline
1 & 3 & 3 & 3 \\
2 & 2 & 2 & 2 \\
3 & 1 & 1 & 1 
\end{array}
\quad
\begin{array}{r|rrr}
\otr & 1 & 2 & 3 \\ \hline
1 & 3 & 1 & 3 \\
2 & 2 & 2 & 2 \\
3 & 1 & 3 & 1 
\end{array}
\quad
\begin{array}{r|c}
x & v(x) \\ \hline
1 & 3 \\
2 & 2 \\
3 & 1 \\
\end{array}.
\]  
Then the tables
\[
\begin{array}{r|rrr}
\phi & 1 & 2 & 3 \\ \hline
1 & 0 & 0 & 2 \\
2 & 6 & 0 & 6 \\
3 & 5 & 5 & 0
\end{array}
\quad
\begin{array}{r|rrr}
\psi & 1 & 2 & 3 \\ \hline
1 & 0 & 1 & 2 \\
2 & 6 & 0 & 1 \\
3 & 5 & 6 & 0
\end{array}
\]
specify a Boltzmann weight with coefficients in $\mathbb{Z}_7$.
The virtual Hopf link then has homset
\[\scalebox{0.85}{\includegraphics{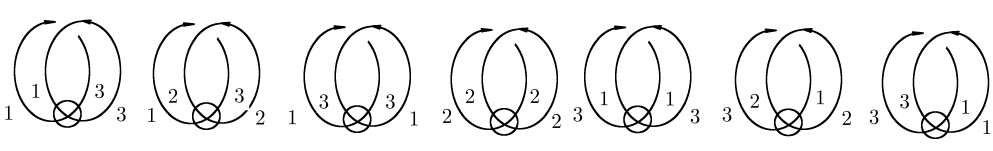}}\]
with Boltzmann weight multiset $\{0,6,0,0,0,6,0\}$.
\end{example}

We can think of Boltzmann weights more algebraically in the following way.
Let $X'=(X\times X)\setminus \Delta(X)$ be the set of ordered pairs of
elements of $X$ minus the diagonal and let $\mathbb{Z}_+$ be the nonnegative
integers. A virtual biquandle colored oriented virtual link diagram determines 
an element of $\vec{v}_c\oplus\vec{v}_v\in 
\mathbb{Z}[X']\oplus\mathbb{Z_+}[X']$ where $\vec{v}_C$ represents the colors of the classical crossings and $\vec{v}_v$ represents the colors of the virtual 
crossings. A Boltzmann weight determines a linear map 
$\phi\oplus \psi:\mathbb{Z}[X']\oplus\mathbb{Z_+}[X']\to A$ in the generic case
or $\phi\oplus \psi:\mathbb{Z}[X']\oplus\mathbb{Z_+}[X']\to A\oplus A$ in the 
strongly compatible case.

\begin{example}
Consider the virtual biquandle $(X,V)$ 
\[
\begin{array}{r|rrrr}
\utr & 1 & 2 & 3 & 4 \\ \hline
1 & 2 & 2 & 2 & 2 \\
2 & 3 & 3 & 3 & 3 \\
3 & 4 & 4 & 4 & 4 \\
4 & 1 & 1 & 1 & 1
\end{array}
\quad
\begin{array}{r|rrrr}
\utr & 1 & 2 & 3 & 4 \\ \hline
1 & 2 & 2 & 2 & 2 \\
2 & 3 & 3 & 3 & 3 \\
3 & 4 & 4 & 4 & 4 \\
4 & 1 & 1 & 1 & 1
\end{array}
\quad
\begin{array}{r|r}
x & V(x) \\ \hline
1 & 3 \\
2 & 4 \\
3 & 1 \\
4 & 2 
\end{array}
\quad
\]
and Boltzmann weight with $A=\mathbb{Z}_5$ coefficients
\[
\begin{array}{r|rrrr}
\phi & 1 & 2 & 3 & 4 \\ \hline
1 & 0 & 0 & 4 & 4 \\
2 & 4 & 0 & 1 & 4 \\
3 & 3 & 3 & 0 & 0 \\
4 & 0 & 3 & 4 & 0
\end{array}
\quad
\begin{array}{r|rrrr}
\psi & 1 & 2 & 3 & 4 \\ \hline
1 & 0 & 3 & 1 & 3 \\
2 & 2 & 0 & 3 & 0 \\
3 & 4 & 2 & 0 & 4 \\
4 & 2 & 0 & 3 & 0
\end{array}.
\]
Then \[X'=\{(1,2),(1,3),(1,4),(2,1),(2,3),(2,4),(3,1),(3,2),(3,4),(4,1),(4,2),(4,3)\}\] and
the $(X,V)$-coloring of the virtual knot $3.1$ 
\[\includegraphics{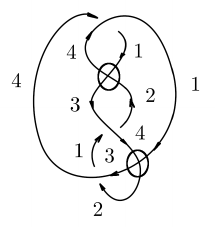}\] is represented by
\[
\left[\begin{array}{cccccccccccc} 0 & 0 & 1 & 0 & 0 & 0 & -1 & 0 & 0 & 0 & 0 & 0\end{array}\right]^T\oplus
\left[\begin{array}{cccccccccccc} 0 & 0 & 2 & 0 & 0 & 0 & 0 & 0 & 0 & 0 & 0 & 0\end{array}\right]^T.
\]
In covector form we can write $\phi\oplus\psi$ as
\[
\left[\begin{array}{cccccccccccc} 0 & 4 & 4 & 4 & 1 & 4 & 3 & 3 & 0 & 3 & 4 & 0\end{array}\right]\oplus
\left[\begin{array}{cccccccccccc} 3 & 1 & 3 & 2 & 3 & 0 & 4 & 2 & 2 & 2 & 0 & 3\end{array}\right].
\]
Then evaluating the Boltzmann weight in the generic case, we compute 
$\phi(\vec{v}_C)+\psi(\vec{v}_V)=(4-3)+(2(3))=1+1=2$; since this weight 
is strongly compatible, we can also keep the classical and virtual parts 
separate, obtaining Boltzmann weight $1\oplus 1$.
\end{example}

\section{\textbf{Quivers and Representations}}\label{QQ}

In this section, we define our new invariants.


\begin{definition}
Let $(X,v)$ be a virtual biquandle, $A$ an abelian group, 
$C=\{(\phi_j,\psi_j)\ |\ j=1,\dots, n\}$ a finite set of virtual
Boltzmann weights, $k$ a choice of coefficient ring and 
$S\subset\mathrm{End}(X,v)$ a set of endomorphisms of $(X,v)$, and let $L$ 
be an oriented classical or virtual link represented by
a virtual link diagram $D$. Then the \textit{virtual biquandle Boltzmann
quiver representation} associated to the \textit{data vector} 
$((X,v),A,C,k,S)$ is the quiver with:
\begin{itemize}
\item A vertex for each virtual biquandle homset element $\vec{v}$,
represented by an $(X,v)$-coloring $D_{\vec{v}}$ of $D$,
\item A copy of the $k$-module $k[A]$ at each vertex (or $k[A]\oplus k[A]$ 
at each vertex if $(\phi_j,\psi_j)$ is strongly compatible for every $j$), 
with distinguished submodule, denoted by $V_C(\vec{v})$, generated by  
$\lbrace \phi_j(\vec{v}) + \psi_j(\vec{v}) \, \vert \, (\phi_j,\psi_j) \in C 
\rbrace$ (or $\lbrace (\phi_j(\vec{v}), \psi_j(\vec{v})) \, \vert \, 
(\phi_j,\psi_j) \in C \rbrace$ if $(\phi_j,\psi_j)$ each pair is strongly 
compatible for every $j$).
\item For each $(\phi,\psi) \in C$, we define a linear transformation 
$f_\sigma: k[A] \rightarrow k[A]$ (or $f_\sigma: k[A] \oplus k[A] \rightarrow 
k[A] \oplus k[A]$ if $(\phi, \psi)$ is strongly compatible) by first setting 
for each $(\phi,\psi) \in C$
\[ f_{\sigma, (\phi,\psi)} (\phi(\vec{v})+\psi(\vec{v})) 
= \phi (\sigma(\vec{v})) + \psi(\sigma(\vec{v})) \]
or 
\[ f_{\sigma, (\phi,\psi)} ( (\phi(\vec{v}),\psi(\vec{v}))) 
= (\phi (\sigma(\vec{v})), \psi(\sigma(\vec{v}))) \]
if $(\phi,\psi) \in C$ is strongly compatible. We then extend linearly,
\[ f_\sigma = \sum_{(\phi,\psi) \in C} f_{\sigma, (\phi,\psi)}.\]
\end{itemize}
\end{definition}

This means that each endomorphism $\sigma \in S$ determines a linear 
transformation $f_\sigma: k[A] \rightarrow k[A]$ sending a distinguished 
submodule to a distinguished submodule.

\begin{definition}
Let $(X,v)$ be a virtual biquandle, $A$ an abelian group, 
$C=\{(\phi_j,\psi_j)\ |\ j=1,\dots, n\}$ a finite set of virtual
Boltzmann weights, $k$ a choice of coefficient ring and 
$S\subset\mathrm{End}(X,v)$ a set of endomorphisms of $(X,v)$. 
Let $L$ be an oriented classical or virtual link. Then we define the 
\textit{virtual biquandle Boltzmann quiver representation of $L$} to be 
the biquandle coloring quiver of $K$ with respect to $((X,v),S)$ with each 
vertex $\vec{v}$ weighted with the pair $(k[A], V_C(\vec{v}))$ and each 
edge defined by $\sigma \in S$ weighted with $f_\sigma$. In the case that 
$(\phi_j, \psi_j)$ is strongly compatible, then the \textit{virtual biquandle 
Boltzmann quiver representation of $K$} to be the biquandle coloring quiver 
of $K$ with respect to $((X,v),S)$ with each vertex $\vec{v}$ weighted with
the pair $(k[A] \oplus k[A], V_C(\vec{v}))$ and each edge 
defined by $\sigma \in S$ weighted with $f_\sigma$.
\end{definition}

By construction, we have our main theorem:

\begin{proposition}
Let $(X,v)$ be a virtual biquandle, $A$ an abelian group, 
$C=\{(\phi_j,\psi_j)\ |\ j=1,\dots, n\}$ a finite set of virtual
Boltzmann weights, $k$ a choice of coefficient ring and 
$S\subset\mathrm{End}(X,v)$ a set of endomorphisms of $(X,v)$. 
Then if two oriented classical or virtual knots or links $L$ and $L'$
are virtually isotopic, their virtual biquandle Boltzmann quiver
representations with respect to the data vector $((X,v),A,C,k,S)$
are isomorphic.
\end{proposition}

\begin{example}\label{Ex:(3.1) and (3.2)}
Let $((X,v),A,C,k,S)$ be the data vector consisting of the biquandle $X$
with operation tables
\[
\begin{array}{r|rrr}
\utr & 1 & 2 & 3 \\ \hline
1 & 1 & 1 & 1 \\
2 & 3 & 2 & 2 \\
3 & 2 & 3 & 3 \\
\end{array}
\quad
\begin{array}{r|rrr}
\otr & 1 & 2 & 3 \\ \hline
1 & 1 & 1 & 1 \\
2 & 3 & 2 & 2 \\
3 & 2 & 3 & 3 \\
\end{array}\quad
\begin{array}{r|c}
x & v(x) \\ \hline
1 & 1 \\
2 & 3 \\
3 & 2 
\end{array}
\]
$A=\mathbb{Z}_3$, 
\[C=\{
[0\ 1\ 1\ 1\ 0\ 1]^T\oplus[0\ 1\ 0\ 1\ 2\ 2]^T,
[0\ 1\ 2\ 0\ 1\ 0]^T\oplus[2\ 0\ 1\ 1\ 0\ 2]^T,
[2\ 2\ 2\ 2\ 2\ 2]^T\oplus[2\ 2\ 1\ 0\ 1\ 0]^T\},\]
$k=\mathbb{Z}$ and $S=\{[1\ 3\ 2],[1\ 1\ 1],[1\ 2\ 3]\}$.
Then we compute the virtual biquandle cohomology quiver representations for
the virtual knots 3.1 and 3.2 respectively as
\[\scalebox{1.2}{\includegraphics{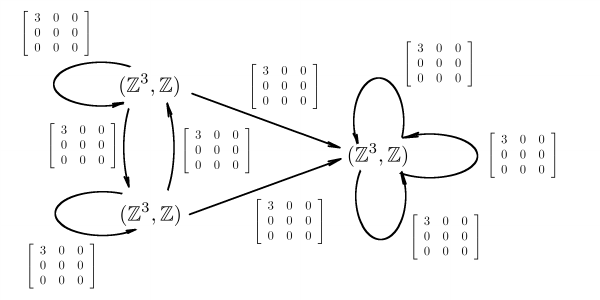}}\]
and
\[\scalebox{1.2}{\includegraphics{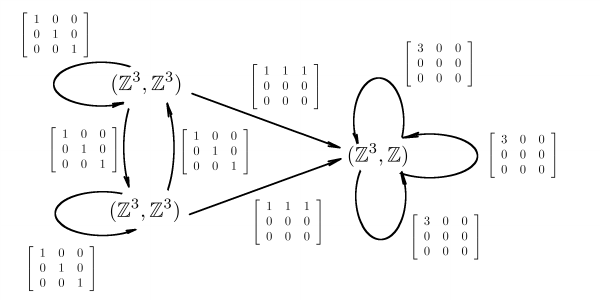}}.\]
\end{example}

Comparing quiver representations directly can get complicated, so as an 
application we define several new infinite families of polynomial invariants
of oriented classical and virtual knots and links from these quivers
via decategorification.

\begin{definition}
Let $L$ be an oriented classical or virtual link and let
$(X,A,C,k,S)$ be a data vector for a virtual biquandle Boltzmann quiver
representation. We define the \textit{virtual biquandle Boltzmann quiver 
representation edge characteristic polynomial} or just the 
\textit{edge characteristic polynomial} to be the sum
\[\Phi_{\chi}^E(L)=\sum_{e\in Q} \chi(M_{\sigma})\]
of the characteristic polynomials of the matrices $M_{\sigma}(x)$
associated to each edge. We further define the
\textit{virtual biquandle Boltzmann quiver representation edge matrix
polynomial} or just the \textit{edge matrix polynomial} to be the sum
\[\Phi_{p_M}^{E}(L)=\sum_{e\in Q} p_M(M_{\sigma})\]
of the matrix polynomials 
\[p_M(M)=\sum_{m_{jk}\in M} m_{jk}x^jy^k\]
(where $m_{jk}$ is the entry in row $j$ column $k$ of M)
encoding the matrices representing $M_{\sigma}$ 
on each edge. In the case of strong compatibility, we can treat the 
classical and virtual matrices separately to get a pair of polynomials
rather than a single polynomial.
\end{definition}

\begin{definition}
Let $L$ be an oriented classical or virtual link and let
$(X,A,C,k,S)$ be a data vector for a virtual biquandle Boltzmann quiver
representation. Then for each maximal directed path $P_j=e_1e_2\dots e_k$ in the
virtual biquandle Boltzmann quiver representation of $L$
with respect to the data vector $((X,v),A,C,k,S)$ without repeated edges, 
in the corresponding sequence of vector spaces and maps
\[V_1\stackrel{f_1}{\to}V_2\stackrel{f_2}{\to}\dots\stackrel{f_{k}}{\to} V_{k+1}\]
let $M_j=f_{k}\dots f_2f_1$ be the resulting matrix product. Then we define the
\textit{virtual biquandle Boltzmann quiver representation maximal path 
characteristic polynomial} or just the \textit{path characteristic polynomial}
to be the sum 
\[\Phi_{\chi}^{P}(L)=\sum_{j} \chi(M_j)s^{|j|}\]
of characteristic polynomials of $M_j$ times a variable 
$s$ to the power of the length of the path $j$ over the set of maximal
nonrepeating paths. We further define the
\textit{virtual biquandle Boltzmann quiver representation maximal path 
matrix polynomial} or just the \textit{path matrix polynomial}
to be the sum
\[\Phi_{p_M}^P(L)=\sum_j p_M(M_j)z^{|j|}\] of products of $p_M$ values of 
$M_j$ times a variable $z$ to the power of the length of the paths $j$ 
over the set of maximal nonrepeating paths in the quiver.
As in the edge polynomial case, in the case of strong compatibility, we can 
treat the classical and virtual matrices separately to obtain a pair of
polynomials in each case.
\end{definition}

We then have, by construction:

\begin{proposition}
The polynomials $\Phi_{\chi}^E,\ \Phi_{\chi}^P,\Phi_{p_M}^E$ and $\Phi_{p_M}^P$ 
are invariants of classical and virtual oriented knots
and links for every choice of finite virtual biquandle $(X,v)$, 
abelian group $A$, set of virtual Boltzmann weights $C$, coefficient ring $k$ 
and set of virtual biquandle endomorphisms $S\subset \mathrm{Hom}(X,X)$.
\end{proposition}

\section{\textbf{Examples and Computations}}\label{EC}

\begin{example}
Let $((X,v),A,C,k,S)$ be the data vector from Example~\ref{Ex:(3.1) and (3.2)}, and we will consider the same two virtual knots from that example, which are $3.1$ and $3.2$. The virtual knot $3.1$ has (via \texttt{python} computation) edge characteristic, edge matrix, path characteristic, path matrix polynomials

\begin{eqnarray*}
    \Phi_\mathcal{X}^E(3.1) &=&9t^3 - 27t^2\\
    \Phi_{p_M}^E(3.1) &=& 27\\
    \Phi_\mathcal{X}^P(3.1) &=& 24s^8t^3 - 157464s^8t^2\\
    \Phi_{p_M}^P(3.1) &=&157464z^8.
\end{eqnarray*}
On the other hand, the virtual knot $3.2$ has the following edge characteristic, edge matrix, path characteristic, path matrix polynomials
\begin{eqnarray*}
    \Phi_\mathcal{X}^E(3.2) &=&9t^3 - 23t^2 + 12t - 4\\
    \Phi_{p_M}^E(3.2) &=& 4x^2y^2 + 4xy + 2y^2 + 2y + 15\\
    \Phi_\mathcal{X}^P(3.2) &=& 24s^8t^3 - 648s^8t^2\\
    \Phi_{p_M}^P(3.2) &=&24z^8(27x^2 + 27x + 27).
\end{eqnarray*}
We can see that all four of the polynomial invariants are sensitive enough to distinguish the virtual knots 3.1 and 3.2.
\end{example}

\begin{example} 

Using \texttt{python} code, we can compute the edge matrix polynomials for the virtual knots with up to 4 classical crossings with a choice of orientation in the table \cite{VKA}. Let $((X,v),A,C,k,S)$ be the data vector consisting of the biquandle $X$
with operation tables
\[
\begin{array}{r|rrrr}
\utr & 1 & 2 & 3 & 4 \\ \hline
1 & 3 & 3 & 3 & 3 \\
2 & 4 & 4 & 4 & 4 \\
3 & 1 & 1 & 1 & 1\\
4 & 2 & 2 & 2 & 2
\end{array}
\quad
\begin{array}{r|rrrr}
\otr & 1 & 2 & 3 & 4 \\ \hline
1 & 3 & 1 & 2 & 4  \\
2 & 2 & 4 & 3 & 1  \\
3 & 4 & 2 & 1 & 3  \\
4 & 1 & 3 & 4 & 2
\end{array}\quad
\begin{array}{r|c}
x & v(x) \\ \hline
1 & 4 \\
2 & 3 \\
3 & 2 \\
4 & 1
\end{array}
\]
$A=\mathbb{Z}_4$, 
\[C=\{
[2\ 3\ 3\ 0\ 3\ 1\ 1\ 3\ 0\ 3\ 3\ 2]^T\oplus[0\ 0\ 0\ 0\ 0\ 0\ 0\ 0\ 0\ 0\ 0\ 0]^T,
[3\ 3\ 2\ 1\ 2\ 1\ 3\ 0\ 1\ 0\ 1\ 3]^T\oplus[2\ 0\ 2\ 2\ 2\ 0\ 0\ 2\ 2\ 2\ 0\ 2]^T\},\]
$k=\mathbb{Z}$ and $S=\{[4\ 3\ 2\ 1], [2\ 1\ 4\ 3], [3\ 4\ 1\ 2] \}$. In this case, both cocycles in $C$ are strongly compatible, so we will compute the edge matrix polynomial for the strongly compatible case. We collect our results in the following table.

\[\begin{array}{r|l}
\Phi_{p_M}^E(L) & \mathrm{Virtual\ knot}\ L \\
\hline
(0, 0) &
2.1, 3.1, 3.2, 3.3, 3.4,
4.2, 4.3, 4.4, 4.5, 4.6, 4.7,
4.10, 4.12, 4.13, 4.14, 4.15, 4.18,\\ &
4.20, 4.22, 4.25, 4.26, 4.27, 4.28,
4.29, 4.30, 4.31, 4.32, 4.33, 4.34,
4.37, 4.38, 4.39, \\ & 4.40, 4.44, 4.46, 
4.48, 4.49, 4.50, 4.51, 4.52, 4.53,
4.54, 4.60, 4.62, 4.63, 4.69, 4.70, \\ &
4.71, 4.73, 4.74, 4.75, 4.78, 4.80, 
4.81, 4.82, 4.83, 4.84, 4.87, 4.88,
4.91, 4.93, \\ &4.94, 4.97, 4.100, 4.101,
4.102, 4.103, 4.104
\\\\
(24, 24) &
3.5, 3.7,
4.8, 4.16, 4.17, 4.19,
4.21, 4.24, 4.41, 4.42,
4.43, 4.45, 4.47, 4.55, \\ &
4.56, 4.57, 4.58, 4.64,
4.66, 4.67, 4.68, 4.72,
4.76, 4.77, 4.79, 4.85, \\ &
4.86, 4.89, 4.90, 4.96,
4.99, 4.105, 4.106, 4.107
\\\\
(12x^{2}y^{2} + 12, 24) &
4.1, 4.9, 4.61, 4.92
\\\\
(12x^{2}y^{2} + 12, 12x^{2}y^{2} + 12) &
4.11, 4.23, 4.35, 4.36, 4.59
\\\\
(36x^{2}y^{2} + 60, 96) &
3.6, 4.108
\\\\
(48x^{2}y^{2} + 48, 24x^{2}y^{2} + 72) &
4.65
\\\\
(24x^{2}y^{2} + 72, 24x^{2}y^{2} + 72) &
4.98
\\\\
(24x^{2}y^{2}, 12x^{2}y^{2} + 12) &
4.95
\end{array}\]

\end{example}

The following example demonstrates that these invariants detect different features of a virtual knot. In particular, virtual knots can agree with respect to one invariant but be distinguished by another. This implies that the invariants are complementary rather than equivalent.

\begin{example}
Let $((X,v),A,C,k,S)$ be the data vector consisting of the biquandle $X$
with operation tables
\[
\begin{array}{r|rrr}
\utr & 1 & 2 & 3 \\ \hline
1 & 1 & 1 & 1 \\
2 & 3 & 2 & 2 \\
3 & 2 & 3 & 3 \\
\end{array}
\quad
\begin{array}{r|rrr}
\otr & 1 & 2 & 3 \\ \hline
1 & 1 & 1 & 1 \\
2 & 3 & 2 & 2 \\
3 & 2 & 3 & 3 \\
\end{array}\quad
\begin{array}{r|c}
x & v(x) \\ \hline
1 & 1 \\
2 & 3 \\
3 & 2 
\end{array}
\]
$A=\mathbb{Z}_3$, 
\[C=\{[2\ 2\ 2\ 2\ 2\ 2]^T\oplus[2\ 2\ 1\ 0\ 1\ 0]^T\},\]
$k=\mathbb{Z}$ and $S=\{[1\ 3\ 2],[1\ 1\ 1],[1\ 2\ 3]\}$. Note that the pair of 2-cocycles is strongly compatible; this means that each invariant in this example will be a pair. We will compute the path characteristic and path matrix polynomials via \texttt{python} code of the virtual knots 3.4 and 4.1 with a choice of orientation. For the virtual knot 3.4, we obtain:
\begin{eqnarray*}
    \Phi_\mathcal{X}^P(3.4) &=& (24 s^{8} t_{1}^{3}, 24 s^{8} (t_{2}^{3} - t_{2}^{2}))\\
    \Phi_{p_M}^P(3.4) &=& (24 x^{2} z^{8}, 24 z^{8}).
\end{eqnarray*}
For the virtual knot 4.1, we obtain
\begin{eqnarray*}
    \Phi_\mathcal{X}^P(4.1) &=& (24 s^{8} t_{1}^{3}, 24 s^{8} (t_{2}^{3} - t_{2}^{2}))\\
    \Phi_{p_M}^P(4.1) &=& (24 x z^{8}, 24 z^{8}).
\end{eqnarray*}

We can see that in this case $\Phi_\mathcal{X}^P(3.4) = \Phi_\mathcal{X}^P(4.1)$, but the max path matrix polynomial can distinguish the pair since $\Phi_{p_M}^P(3.4) \neq \Phi_{p_M}^P(4.1)$.
\end{example}

\section{\textbf{Questions}}\label{Q}

We end with some questions and directions for future research. 

It is important to note that the examples in this paper are toy examples
meant to illustrate the computation of the invariant using small
biquandles and finite rings easily accessible by python computation.
The true power of this infinite family of invariants lies in using
larger biquandles and larger finite or infinite rings. In particular,
fast algorithms for computing these invariants is of great interest.

What is the geometric meaning of these new polynomial invariants? What 
properties characterize the polynomial values? What are the relationships
with other families of virtual knot and link invariants?

\bibliography{ab-jc-sn}{}

@misc{VKA,
    AUTHOR = {Green, Jeremy},
     TITLE = {A Table of Virtual Knots \textup{https://www.math.toronto.edu/drorbn/Students/GreenJ/}},
   JOURNAL = {},
  FJOURNAL = {},
    VOLUME = {},
      YEAR = {},
    NUMBER = {},
     PAGES = {},
      ISSN = {},
     CODEN = {},
   MRCLASS = {},
  MRNUMBER = {},
MRREVIEWER = {},
       DOI = {},
       URL ={https://www.math.toronto.edu/drorbn/Students/GreenJ/},
}

@article {JCN1,
    AUTHOR = {Ceniceros, Jose and Nelson, Sam},
     TITLE = {Virtual {Y}ang-{B}axter cocycle invariants},
   JOURNAL = {Trans. Amer. Math. Soc.},
  FJOURNAL = {Transactions of the American Mathematical Society},
    VOLUME = {361},
      YEAR = {2009},
    NUMBER = {10},
     PAGES = {5263--5283},
      ISSN = {0002-9947,1088-6850},
   MRCLASS = {57M27 (18G60)},
  MRNUMBER = {2515811},
       DOI = {10.1090/S0002-9947-09-04751-5},
       URL = {https://doi-org.ccl.idm.oclc.org/10.1090/S0002-9947-09-04751-5},
}

@article{CN,
author = {Cho, Karina and Nelson, Sam},
title = {Quandle coloring quivers},
journal = {Journal of Knot Theory and Its Ramifications},
volume = {28},
number = {01},
pages = {1950001},
year = {2019},
doi = {10.1142/S0218216519500019},
URL = { https://doi.org/10.1142/S0218216519500019},
eprint = {https://doi.org/10.1142/S0218216519500019},
}

@article {CN2,
    AUTHOR = {Cho, Karina and Nelson, Sam},
     TITLE = {Quandle cocycle quivers},
   JOURNAL = {Topology Appl.},
  FJOURNAL = {Topology and its Applications},
    VOLUME = {268},
      YEAR = {2019},
     PAGES = {106908, 10},
      ISSN = {0166-8641,1879-3207},
   MRCLASS = {57M27 (57M25)},
  MRNUMBER = {4018585},
MRREVIEWER = {Markus\ Szymik},
       DOI = {10.1016/j.topol.2019.106908},
       URL = {https://doi.org/10.1016/j.topol.2019.106908},
}

@book {EN,
    AUTHOR = {Elhamdadi, Mohamed and Nelson, Sam},
     TITLE = {Quandles---an introduction to the algebra of knots},
    SERIES = {Student Mathematical Library},
    VOLUME = {74},
 PUBLISHER = {American Mathematical Society, Providence, RI},
      YEAR = {2015},
     PAGES = {x+245},
      ISBN = {978-1-4704-2213-4},
   MRCLASS = {57M27 (57M25 57Q45)},
  MRNUMBER = {3379534},
}

@article {J,
    AUTHOR = {Joyce, David},
     TITLE = {A classifying invariant of knots, the knot quandle},
   JOURNAL = {J. Pure Appl. Algebra},
  FJOURNAL = {Journal of Pure and Applied Algebra},
    VOLUME = {23},
      YEAR = {1982},
    NUMBER = {1},
     PAGES = {37--65},
      ISSN = {0022-4049},
     CODEN = {JPAAA},
   MRCLASS = {57M25 (20F29 20N05 53C35)},
  MRNUMBER = {638121 (83m:57007)},
MRREVIEWER = {Mark E. Kidwell},
       DOI = {10.1016/0022-4049(82)90077-9},
       URL = {http://dx.doi.org/10.1016/0022-4049(82)90077-9},
}

@article {K,
    AUTHOR = {Kauffman, Louis H.},
     TITLE = {Virtual knot theory},
   JOURNAL = {European J. Combin.},
  FJOURNAL = {European Journal of Combinatorics},
    VOLUME = {20},
      YEAR = {1999},
    NUMBER = {7},
     PAGES = {663--690},
      ISSN = {0195-6698},
   MRCLASS = {57M25 (57M27)},
  MRNUMBER = {1721925 (2000i:57011)},
MRREVIEWER = {Olivier Collin},
       DOI = {10.1006/eujc.1999.0314},
       URL = {http://dx.doi.org/10.1006/eujc.1999.0314},
}

@article {KM,
    AUTHOR = {Kauffman, Louis H. and Manturov, Vassily O.},
     TITLE = {Virtual biquandles},
   JOURNAL = {Fund. Math.},
  FJOURNAL = {Fundamenta Mathematicae},
    VOLUME = {188},
      YEAR = {2005},
     PAGES = {103--146},
      ISSN = {0016-2736,1730-6329},
   MRCLASS = {57M25 (57M27)},
  MRNUMBER = {2191942},
MRREVIEWER = {Heather\ A.\ Dye},
       DOI = {10.4064/fm188-0-6},
       URL = {https://doi-org.ccl.idm.oclc.org/10.4064/fm188-0-6},
}

@misc{SN1,
    AUTHOR = {Nelson, Sam},
     TITLE = {Quandle Cohomology Quiver Representations},
   JOURNAL = {Preprint. arXiv:2411.02153},
  FJOURNAL = {arXiv:2411.02153},
    VOLUME = { },
      YEAR = {Preprint. arXiv:2411.02153, 2024},
    NUMBER = { },
     PAGES = { },
      ISSN = {},
     CODEN = {},
   MRCLASS = {},
  MRNUMBER = {},
MRREVIEWER = {},
       DOI = {},
       URL ={https://arxiv.org/abs/2411.02153},
}

@article {NT,
    AUTHOR = {Nelson, Sam and Tamagawa, Sherilyn},
     TITLE = {Quotient quandles and the fundamental {L}atin {A}lexander
              quandle},
   JOURNAL = {New York J. Math.},
  FJOURNAL = {New York Journal of Mathematics},
    VOLUME = {22},
      YEAR = {2016},
     PAGES = {251--263},
      ISSN = {1076-9803},
   MRCLASS = {57M27 (57M25)},
  MRNUMBER = {3484684},
MRREVIEWER = {Zhiyun\ Cheng},
       URL = {http://nyjm.albany.edu:8000/j/2016/22_251.html},
}
\bibliographystyle{abbrv}

\bigskip

\noindent
\textsc{Department of Mathematical Sciences \\
Claremont McKenna College \\
850 Columbia Ave. \\
Claremont, CA 91711} 

\medskip

\noindent
\textsc{Mathematics and Statistics Department \\
Hamilton College \\
198 College Hill Rd. \\
Clinton, NY 13323}

\end{document}